\documentclass[a4paper,10pt]{article}
\usepackage[french]{babel}
\usepackage{subfigure}
\usepackage{enumerate,amsmath,amsthm,amssymb,mathrsfs}
\usepackage{graphicx,color,upgreek}
\usepackage[latin1]{inputenc}
 \usepackage[T1]{fontenc}

\def\a{\alpha}

\def\b{\beta}

\def\g{\gamma}

\def\at{\alpha^\vee}
\def\bt{\beta^\vee}
\def\bti{\beta_i^\vee}
\def\ati{\alpha_i^\vee}
\def\gt{\gamma^\vee}

\def\lieg{\mathfrak g}
\def\lieh{\mathfrak h}
\def\car{{\rm car\,}}
\newcommand\scal[2]{\langle #1 , #2\rangle}

\def\C{{\mathbb C}}

\def\Z{{\mathbb Z}}
\def\R{{\mathbb R}}
\def\N{{\mathbb N}}
\def\Q{{\mathbb Q}}

\newenvironment{preuve}{{\noindent\sc Preuve.}}{\hfill $\square$\vspace*{.5cm}}

\theoremstyle{plain}
\newtheorem{Th}{Th\'eor\`eme}[section]

\newtheorem{Prop}[Th]{Proposition}
\newtheorem{Cor}[Th]{Corollaire}

\theoremstyle{definition}
\newtheorem{Df}{D\'efinition}[section]
\newtheorem{Ex}{Exemple}
\newtheorem{Exs}[Ex]{Exemples}

\theoremstyle{remark}
\newtheorem*{NB}{Remarque}
\newtheorem*{NBs}{Remarques}

\date{}

\begin{document}
\title{Composantes PRV généralisées et chemins de Littelmann}
\author{Pierre-Louis Montagard\footnote{Université Montpellier II - 
CC 51-Place Eugène Bataillon -
34095 Montpellier Cedex 5 -
France - {\tt pierre-louis.montagard@math.univ-montp2.fr}}}

\maketitle
{\small 
\begin{center}
{\bf Abstract}

\end{center}

We give a sufficient condition for a Littelmann path
to represent a vector of extremal weight of an integrable irreducible highest
weight representation of a symmetrisable Kac-Moody
algebra. Thanks to this condition we present, in a more general
context, 
an alternative proof of a recent result by Boris
Pasquier, Nicolas Ressayre and 
the author of this article on the existence of generalized PRV components. 

\begin{center}
{\bf Résumé}
\end{center}
Nous énonçons une condition suffisante pour qu'un
chemin de Littelmann représente un vecteur de poids extrémal d'une
représentation intégrable, irréductible et de plus haut poids d'une
algèbre de Kac-Moody symétrisable. \`A 
l'aide de 
cette condition, nous présentons, dans un contexte plus général, une
preuve alternative de résultats de Boris
Pasquier, Nicolas Ressayre et l'auteur de cet article sur l'existence
de  composantes PRV généralisées.}

\section{Introduction}
La conjecture PRV a été énoncée dans les années 60 dans
\cite{PRV:original} par  Parthasarathy,
Ranga Rao et Varadarajan. Cette conjecture concerne le
problème de décomposition du produit tensoriel de deux représentations
irréductibles d'une algèbre semi-simple complexe $\lieg$. Son énoncé
est très simple~: soient $\mu$ et
$\nu$ deux poids dominants de $\lieg$ et $v,w$ deux
éléments du groupe de Weyl de $\lieg$ ; soient $V(\mu)$
(resp. $V(\nu)$) la représentation
irréductible de plus haut poids $\mu$
(resp. $\nu$). Alors si le poids $\lambda=v\mu+w\nu$ est dominant la
représentation irréductible de plus haut poids $\lambda=v\mu+w\nu$ est de
multiplicité non nulle dans
$V(\mu)\otimes V(\nu)$. 

On connaît aujourd'hui plusieurs preuves de cette conjecture. Citons
d'abord, les preuves (simultanées et indépendantes) à la fin des
années 80, de Shrawan
Kumar \cite{Kumar:prv1} et \cite{Kumar:prv2} et Olivier
Mathieu \cite{Mathieu:prv}. La géométrie, et notamment la géométrie des
espaces de drapeaux du 
groupe $G$ associé à $\lieg$ intervient de façon essentielle dans ces
deux preuves. Notons que ces deux preuves sont valables dans le
contexte plus général des algèbres de Kac-Moody symétrisables.
Plus tard, au milieu des années 90, Peter Littelmann a donné une preuve
combinatoire et plus élémentaire de la conjecture PRV. Cette preuve est une
application de la théorie des chemins de Littelmann développée dans
\cite{pet1} et \cite{pet2}, théorie qui généralise la
règle de Littlewood-Richardson dans le cadre des algèbres de Kac-Moody
symétrisables.  

Récemment dans \cite{MPR1} et \cite{MPR2}, avec Nicolas Ressayre et 
Boris Pasquier, nous avons démontré plusieurs généralisations de
l'énoncé PRV dans le contexte d'un groupe
algébrique réductif $G$ et en utilisant des outils géométriques.
Énonçons l'une de ces généralisations. 

\begin{Th}\label{th:prv:intro}
  Soit $G$ une groupe algébrique réductif complexe, soient
  $(\mu,\nu)$ un couple de poids dominants et  $(v,w)$ un couple
d'éléments du groupe de Weyl,  soit $\beta$ une racine
positive (et $\bt$ la coracine associée) telle que l'une des
conditions suivantes soit
satisfaite~:
\begin{enumerate}[(i)]
\item $\beta$ est une racine simple ;
\item $v^{-1}\beta$ est une racine simple ; 
\item $w^{-1}\beta$ est une racine simple.
\end{enumerate}
Soit $k$ un entier vérifiant les inégalités~: 
$$0\leq k\leq\min\{\scal{v\mu}\bt,\scal{w\nu}\bt\},$$ 
alors si le poids
$\lambda=v\mu+w\nu-k\beta$ est dominant la représentation
irréductible $V(\lambda)$ est de multiplicité non nulle dans le
produit tensoriel $V(\mu)\otimes V(\nu)$. 
\end{Th}

L'ensemble des composantes ainsi obtenues contenant strictement les
composantes obtenues par l'énoncé de la conjecture PRV originale, nous les
appellerons composantes PRV généralisées.   
Nous avons représenté dans la figure~1, 
un exemple, montrant les composantes obtenues par ce procédé. 
Dans cet exemple $G$ est égal au groupe ${\rm Sl}_3(\C)$,
$\mu=7\varpi_1+3\varpi_2$, $\nu=\varpi_1+3\varpi_2$ où
$\varpi_1,\varpi_2$ sont les deux poids fondamentaux de $G$. Nous avons représenté sur la figure les
poids dominants correspondants aux représentations irréductibles de
multiplicité non nulle du produit tensoriel $V(\mu)\otimes V(\nu)$, les composantes PRV classiques, les
composantes PRV généralisées ainsi que les segments de direction
$\beta$ donnés par le théorème \ref{th:prv:intro} et contenant une
composante PRV généralisée (et non classique).  Le réseau représenté
est le réseau engendré par les racines (et non le réseau des
poids). En effet, comme le
poids $\mu+\nu$ appartient au réseau des racines, toutes
les composantes qui apparaissent dans le produit tensoriel
$V(\mu)\otimes V(\nu)$ appartiennent à ce réseau.

Dans cet article, nous allons donner une preuve de l'existence des
composantes PRV généralisées pour les algèbres de Kac-Moody, en
utilisant la théorie des chemins de Littelmann. 
Nous ne retrouvons pas exactement le même énoncé dans le
contexte des algèbres de Kac-Moody symétrisable. En effet, dans le cas
{\it (i)} ($\beta$
est une racine simple) la preuve du théorème n'est
valable que si $\lieg$ est de dimension finie et
nous ne savons pas si ce cas est vrai pour les algèbres de Kac-Moody
symétrisables.

Rappelons brièvement les principaux résultats de la théorie des
chemins de Littelmann,
résultats qui seront détaillés et précisés dans la
partie \ref{part:1}. 
Si $X$ est le réseau des poids de $\lieg$,
on considère les chemins tracés dans $X\otimes_\Z\R$
d'origine $0$ et d'extrémité un élément de $X$. Si $\mu$ est un
poids dominant, et si $\pi$ est un chemin d'extrémité
$\mu$ complètement contenu dans la chambre dominante, la théorie des chemins
de Littelmann permet de définir un ensemble de chemins
$B(\pi)$ contenant $\pi$ tels que la multiplicité du
poids $\chi$ dans $V(\mu)$ est égale au nombre de chemins dans
$B(\pi)$ d'extrémité $\chi$ (théorème
\ref{th:car:chemin}).

Le modèle des chemins donne également une règle combinatoire pour
décomposer le
produit tensoriel de deux représentations irréductibles. En effet, si
$\nu$ est un poids dominant de
$\lieg$, si $\pi'$ est un chemin contenu dans la
chambre dominante et d'extrémité $\nu$ et si
$B(\pi)*B(\pi')$  est l'ensemble des
concaténations des chemins de $B(\pi)$ et des chemins de
$B(\pi')$, alors la multiplicité
de la représentation irréductible $V(\lambda)$ dans $V(\mu)\otimes
V(\nu)$ est égale au nombre de chemins dans $B(\pi)*B(\pi')$
contenus dans la chambre dominante et 
d'extrémité $\lambda$ (voir le théorème \ref{th:decompo:chemin}). 

Rappelons que pour tout $w\in W$ le
poids $w\lambda$ est de multiplicité un dans $V(\lambda)$ ; pour tout
$w\in W$, il existe
donc un unique chemin dans $B(\pi)$  d'extrémité
$w\lambda$ ; ces chemins seront appelés chemins extrémaux. Le chemin
extrémal correspondant à $w$ égal à l'identité
est caractérisé comme étant l'unique chemin de $B(\pi)$
contenu dans la chambre dominante.
Pour $w$ différent de l'identité, on ne connaît
pas de caractérisation du même genre pour les chemins extrémaux
d'extrémité $w\lambda$. Dans la
partie \ref{part:2}, nous énonçons une condition suffisante pour qu'un chemin
soit extrémal. 

Enfin, dans la partie \ref{part:3}, nous montrons le théorème principal
de ce  travail~: l'existence des composantes PRV généralisées pour les
algèbres de Kac-Moody symétrisables. Pour cela, nous utilisons le
critère obtenu dans la partie  \ref{part:2} pour montrer qu'un chemin
explicite de $B(\pi)*B(\pi')$ est extrémal. 
Nous énonçons et démontrons ensuite un résultat  plus général
(le théorème \ref{th:prvbis}), ou nous exhibons des composantes
du même type que les composantes PRV généralisées, mais 
qui dépendent d'un ensemble de racines deux à deux
orthogonales. Ce résultat est également une généralisation d'un
résultat de l'article \cite{MPR2}. 

\section{Rappels}\label{part:1}
Pour fixer les notations nous allons rappeler brièvement les points
principaux de la théorie
des chemins. Tous ces résultats sont dus à Peter Littelmann et nous renvoyons aux articles originaux de
\cite{pet1} et \cite{pet2} pour les preuves et les
détails.
Dans tout cette article $\lieg$ désignera une algèbre de Kac-Moody
symétrisable sur le corps des nombres complexes. Nous utiliserons les notation suivantes concernant
$\lieg$~:
\begin{itemize}
\item $\lieh$ : une sous-algèbre de Cartan ;
\item $X$ : le réseau des poids ;
\item $X_\Q=X\otimes_\Z\Q$ et $X_\R=X\otimes_\Z\R$ ;
\item $S$ : l'ensemble des racines simples ;
\item $W$ : le groupe de Weyl ;
\item si $\beta$ est une racine réelle, nous noterons
  $\bt\in X_\R^*$ la
  coracine associée et si $\chi\in X_\R$ nous noterons $\scal\chi\bt$
  l'évaluation $\bt(\chi)$ ; 
\item si $\a\in S$, $s_\a$ désigne la réflexion sur $X_\R$ définie par
  $s_\a(\chi)=\chi-\scal\chi\at\a$ ;
\item  $D$ désignera la chambre dominante et si $\lambda\in X\cap D$
  nous noterons $V(\lambda)$ la représentation intégrable irréductible
  et de plus haut poids $\lambda$.
\end{itemize}

\begin{Df}
  Un chemin $\pi$ est une application $\pi\, :\, [0,1]\rightarrow X_\R$
  rectifiable telle que $\pi(0)=0$ et $\pi(1)\in X$. On
dit que deux chemins $\pi$ et $\pi'$ sont équivalents s'il existe une
reparamétrisation $\phi\, :\, [0,1]\rightarrow[0,1]$ croissante,
surjective et continue telle que $\pi=\pi'\circ\phi$. Nous
considérerons les chemins modulo l'équivalence ci-dessus. 
\end{Df}

\begin{Exs}
  \begin{enumerate}
  \item Soit $\chi\in X$, nous noterons $\pi_\chi$ le chemin défini
    par $\pi_\chi(t)=t\chi$.
\item Si $\pi_1$ et $\pi_2$ sont deux chemins, le chemin $\pi_1*\pi_2$
  est le chemin défini par~:
$$
  \begin{array}{rcll}
    \pi(t)&=&\pi_1(2t)&\ {\rm pour}\  t\leq 1/2 \; ;\\
 &=&\pi_1(1)+\pi_2(2t-1)&\ {\rm pour} \ t>1/2.
   \end{array}
$$
\item Dans la suite, nous considérerons essentiellement des chemins affines par morceaux et tels que
  les changements de direction se font en des points rationnels. Un
  tel chemin $\pi$ s'écrit
  $\pi=\pi_{\chi_1}*\pi_{\chi_2}*\cdots*\pi_{\chi_p}$ avec
  $(\chi_1,\chi_2,\ldots,\chi_p)\in X_\Q^p$.   
\item Si $\pi$ est un chemin, nous noterons $\pi^*$ le chemin défini
  par $\pi^*(t)=\pi(1-t)-\pi(1)$. 
\end{enumerate} 
\end{Exs}

Pour toute racine réelle $\alpha$ nous noterons la fonction :
$H_\a^\pi(t)=\scal{\pi(t)}\at$. Soit $m_\alpha^\pi$ le minimum de
cette fonction. Nous aurons également besoin des fonctions suivantes~:
$$
L^\pi_\a(t)=\min\{1,(H_\a^\pi(s)-m_\a^\pi)_{t\leq s\leq 1}\}
$$
et 
$$
R^\pi_\a(t)=\max\{0,(m^\pi_\a-H_\a^\pi(s))_{0\leq s\leq t}\}.
$$
\begin{Df}
  Soit $\pi\in \Pi$ et $\a$ une racine simple, alors si $L_\a^\pi(1)<1$,
  $f_\a\pi$ n'est pas défini, sinon
$f_\a\pi(t):=\pi(t)-L^\pi_\a(t)\a$.

De même si $R^\pi_\a(1)>0$ alors $e_\a\pi$ n'est pas défini, sinon $e_\a\pi(t):=\pi(t)-R^\pi_\a(t)\a$. 
\end{Df}

Voici les propriétés élémentaires des opérateurs $f_\a$ et $e_\a$~:

\begin{Prop}\label{prop:simple}
Soit $\pi\in \Pi$ et $\a$ une racine simple ;
  \begin{enumerate}
  \item Si $f_\a\pi$ est défini alors $f_\a\pi(1)=\pi(1)-\a$. De même
    si $e_\a\pi$ est défini, alors $e_\a\pi(1)=\pi(1)+\a$. 
\item Si $f_\a\pi$ est défini, alors $e_\a f_\a\pi$ est défini et
    $e_\a f_\a\pi=\pi$. De même, si  $e_\a\pi$ est défini, alors $f_\a e_\a\pi$ est défini et
    $f_\a e_\a\pi=\pi$.
\item Si $e_\a\pi$ est défini, alors $f_\a\pi^*$ est défini et
  $f_\a\pi^*=(e_\a\pi)^*$. De même, si $f_\a\pi$ est défini, alors $e_\a\pi^*$ est défini et
  $e_\a\pi^*=(f_\a\pi)^*$.
\item Soit $m$ (resp. $n$), maximal tel que $f_\a^m\pi$
  (resp. $e_\a^n\pi$) soit  défini, alors $n-m=\scal{\pi(1)}\a$,
$$
m=\max\{a\in\Z\,|\,a\leq \scal{\pi(1)}\at -m_\a^\pi\}\  {\rm
  et}\  n=\max\{a\in\Z\,|\,a\leq |m_\a^\pi|\}.
$$

\item  Si $e_\a\pi$ (resp. $f_\a\pi$) est défini alors pour
  tout $n\in\N$, $e^n_\a (n\pi)$ (resp. $f^n_\a (n\pi)$) est défini et on
  a : $e^n_\a (n\pi)=ne_\a\pi$ (resp. $f^n_\a (n\pi)=nf_\a\pi$).

  \end{enumerate}
\end{Prop}

Si $\pi\in\Pi$, nous noterons $B(\pi)$ le plus petit sous-ensemble de
$\Pi$  contenant $\pi$ et stable par les opérateurs $e_\a$ et $f_\a$
(pour $\a\in S$).
Si $B$ est un sous-ensemble de $\Pi$, nous définirons $\car
B=\sum_{\pi\in B}e^{\pi(1)}$ ; nous dirons que $B$ est entier si pour
tout $\pi\in B$ et pour toute racine $\a\in S$, le minimum $m_\a^\pi$
est un entier. Enfin nous dirons qu'un chemin $\pi$ est dominant si
son image est contenue dans la chambre dominante $D$ et nous noterons
$\Pi^+$ l'ensemble des chemins dominants.

Nous pouvons maintenant énoncer les résultats fondamentaux de
P. Littelmann.

\begin{Th}\label{th:car:chemin}
  Soit $\pi\in\Pi$ un chemin dominant, alors on a les résultats suivants :
\begin{itemize}
\item l'ensemble $B(\pi)$ est entier ;
\item $\pi$ est l'unique chemin dominant de $B(\pi)$ ;
\item $\car B(\pi)=\car V(\pi(1))$. 
\end{itemize}
\end{Th}

Si $B$ et $B'$ sont deux sous-ensemble de $\Pi$ nous noterons $B*B'$ l'ensemble des concaténations
$$
B*B':=\{\pi*\pi'\ |\pi\in B,\, \pi'\in B'\}.
$$ 

\begin{Th}\label{th:decompo:chemin}
Soit $\pi$ et $\pi'$ deux chemins dominants, alors l'ensemble
$B(\pi)*B(\pi')$ est entier et se décompose en union disjointe~:
$$B(\pi)*B(\pi')=\bigcup_{\pi*\eta\in\Pi^+,\,\eta\in B(\pi')}B(\pi*\eta).$$
\end{Th}

Enfin, on a la conséquence suivante qui permet de décomposer le produit
tensoriel de deux représentations irréductibles de $\lieg$.

\begin{Th}\label{th:decompo:representation}
  Soit $\lieg$ une algèbre de Kac-Moody, $\mu$ et $\nu$ deux poids
  dominants. Soit $\pi$ un chemin dominant (resp. $\pi'$) tel que
  $\pi(1)=\mu$ (resp. $\pi'(1)=\nu$), alors on a la décomposition
  suivante comme $\lieg$-module~:
$$
V(\mu)\otimes V(\nu)=\bigoplus_{\pi*\eta\in\Pi^+,\,\eta\in
  B(\pi')}V(\mu+\eta(1)).
$$  
\end{Th}

Pour finir cette partie, rappelons que les opérateurs $e_\a$ et $f_\a$
permettent de définir une action du
groupe de Weyl sur $\Pi$. 

\begin{Df}\label{df:weyl}
 Pour tout $\pi\in\Pi$
et pour toute racine simple $\a\in S$, définissons~:
$$
\tilde s_\a\pi:=\left\{
\begin{array}{rl}
  &f^n_\a(\pi);\ {\rm si}\ n=\scal{\pi(1)}\at\geq 0\\
  &e^{-n}_\a(\pi);\ {\rm sinon.}
\end{array}
\right.
$$
\end{Df}

\begin{Th}\label{prop:weyl}
  L'application $s_\a\mapsto \tilde s_\a$ définit une action de $W$ sur
  $\Pi$.
\end{Th}

\section{Chemins extrémaux}\label{part:2}

Rappelons que si $V(\lambda)$ est une représentation irréductible de
$\lieg$, alors pour tout $w\in W$, le poids $w\lambda$ est un poids de
multiplicité un de $V(\lambda)$ appelé poids extrémal de $V(\lambda)$. 
Nous allons maintenant définir la notion de chemin extrémal. 
\begin{Df}
  Soit $\eta$ un chemin, et soit $\pi$ l'unique chemin dominant tel
  que $\eta\in B(\pi)$, on dit que $\eta$ est un chemin extrémal si le
  poids $\eta(1)$ est un poids extrémal de la représentation
  irréductible $V(\pi(1))$. 
\end{Df}

Si $\lambda$ est un poids dominant, et si $\pi=\pi_\lambda$ est le
chemin direct entre $0$ et un poids dominant $\lambda$ alors
l'ensemble des chemins de $B(\pi_\lambda)$ sont appelés chemins de
Lakshmibai-Seshadri et admettent une description combinatoire
explicitée par Littelmann (voir \cite{pet1}). Indépendamment de cette
classification, il est
facile de décrire les chemins extrémaux de $B(\pi_\lambda)$. 

\begin{Prop}\label{prop:ls}
Soit $\lambda$ un poids dominant et $w\in W$ un élément du groupe de
  Weyl, alors le chemin $\pi_{w\lambda}$ appartient à
  $B(\pi_\lambda)$ et l'ensemble 
$$\{\pi_{w\lambda}\ |\ w\in W\}$$
est l'ensemble des chemins extrémaux de $B(\pi_\lambda)$. 
\end{Prop}

\begin{preuve}
  On vérifie que si $\a\in S$ et $w\in W$, alors
  $\tilde s_\a\pi_{w\lambda}=\pi_{s_\a w\lambda}$. Comme l'ensemble
  $\{s_\a\,|\, \a\in S\}$ engendre $W$, la proposition s'en
  déduit immédiatement. 
\end{preuve} 

\begin{NB}
Si $\pi$ est un chemin dominant quelconque on ne sait pas, en général, déterminer
  les chemins extrémaux de $B(\pi)$. Cependant dans le cas $A_n$, si
  $\lambda$ est un poids dominant (et donc une partition), on
  peut définir une injection de l'ensemble $T(\lambda)$ des tableaux
  semi-standard de forme $\lambda$  vers l'ensemble des chemins (voir \cite{pet3}). 
 Dans ce contexte, les éléments de $T(\lambda)$ de poids extrémal sont
 les tableaux clefs qui ont été définis et étudiés par Lascoux et
 Schützenberger, notamment dans \cite{LasSch}.  

\end{NB}

Avant d'énoncer un critère qui assure qu'un chemin
est extrémal, rappelons la notion d'ensemble d'inversion
d'un élément du groupe de Weyl~: soit $R^+_{re}$ l'ensemble des racines
réelles positives, et soit $w\in W$ un élément du groupe de Weyl, on
note $I(w)$, l'ensemble
  d'inversion de $w$ :
$$
I(w)=\{\beta \in R^+_{re}\, | \,w\beta\in -R^+_{re}\}. 
$$
Rappelons également que si $\b$ est une racine réelle et $\pi$ est
un chemin, nous avons
défini
$H^\pi_\b(t)=\scal{\pi(t)}\bt$. 

\begin{Th}\label{th:main}
  Soit $\pi$ un chemin tel que pour toute racine réelle positive $\beta$~:
 \begin{itemize}
    \item ou bien pour tout $t\in [0,1]$ $H^\pi_\beta(t)\geq 0$  ;
\item  ou bien, il existe un réel $t^\pi_\beta\in[0,1[$
    tel que la fonction $H^\pi_\beta$ soit positive ou nulle pour $t\leq t^\pi_\beta$  et
    strictement négative et décroissante pour
    $t>t^\pi_\beta$,
\end{itemize}
alors,
\begin{enumerate}[(i)]
\item $\pi$ est un chemin extrémal ;
\item si $w$ est l'élément de
plus petite longueur tel que $w\pi(1)$ est dominant, alors 
$$
I(w)=\{\beta\in R^+_{re}\ |\ H^\pi_\beta(1)<0\}.
$$ 
\end{enumerate}
\end{Th}

\begin{preuve}
 Soit $\pi$ vérifiant l'hypothèse
 du théorème. 
Il existe un unique $w\in W$ de longueur minimale tel que $w\pi(1)$
soit dominant. Nous allons montrer le théorème par
 récurrence sur la longueur de $w$. 

Si $w$ est de longueur $0$, alors $w$ est égal à l'identité,
$\pi(1)$ est 
un poids dominant et donc par hypothèse, pour tout
racine réelle positive 
  $\b$ et pour tout $t\in[0,1]$, $H_\beta^\pi(t)\geq 0$ ; $\pi$ est
 un chemin dominant, il est donc extrémal et on a bien l'égalité
 $I(w)=\emptyset$. 

Si $l(w)>0$, alors il existe une racine simple $\alpha$ telle
que $H^\pi_\a(1)<0$. On a alors~:
$I(w)=s_\a(I(ws_\a))\cup\{\a\}$, et donc $l(ws_\a)=l(w)-1$. Nous
allons montrer que le chemin $\pi'=\tilde s_\a\pi$ vérifie les
hypothèses du théorème. Ce qui permet de conclure, en effet comme
$l(ws_\alpha)<l(w)$,  $ws_\a$ est de longueur minimale tel que
$ws_\a\pi'(1)$ est dominant ; nous obtenons
le point $(i)$ en utilisant l'hypothèse de récurrence. Pour le point
$(ii)$, il faut remarquer de plus que~:
$$
\{\beta\in R^+_{re}\ |\ H^\pi_\beta(1)< 0\}=s_\a(\{\beta\in
R^+_{re}\ |\ H^{\tilde s_\alpha \pi}_\beta(1)< 0\})\cup\{\a\}. 
$$

Montrons donc que $\pi'=\tilde s_\a\pi$ vérifie les hypothèses du
théorème. Quitte à reparamétriser $\pi$, on peut supposer que 
$t^\pi_\alpha=1/2$. Définissons $\pi_1(t)=\pi(t/2)$ et
$\pi_2(t)=\pi((t+1)/2)-\pi(1/2)$. Par définition on a
$\pi=\pi_1*\pi_2$. Rappelons que si $n=\scal{\pi(1)}{\at}$, alors
d'après la définition \ref{df:weyl}, on a~: $\pi'=\tilde s_\a\pi=e_\a^{-n}\pi$
; on en déduit que
$\tilde s_\alpha\pi=\pi_1*s_\alpha(\pi_2)$. Si $\g$ est une racine
réelle, alors un simple calcul montre que~:
$$
\begin{array}{rcll}\label{formule}
H^{\tilde s_\alpha\pi}_\g(t)&=&H^\pi_\g(t)&  \ {\rm si}\ t\leq 1/2\; ;\\
 &=&H^\pi_{s_\a\g}(t)&\ {\rm si}\  t\geq 1/2.
\end{array}
$$
Remarquons que les deux expressions ci-dessus coïncident en $t=1/2$
puisque $H_\a^\pi(1/2)=0$. 
Soit $\gamma$
une racine réelle positive telle que $H^{\tilde s_\a\pi}_\g(1)\geq 0$,
alors ou bien $\g=\alpha$ et dans ce cas pour $t\leq 1/2$, on a 
$H^{\tilde s_\alpha\pi}_{\alpha}(t)=H^{\pi}_{\alpha}(t)$ qui est positive
et pour $t>1/2$, on a $H^{\tilde s_\alpha\pi}_{\a}(t)=-H^\pi_\a(t)$
qui est également positive. 
Si $\gamma\neq\alpha$, alors il existe une racine réelle positive
$\beta$ telle que $s_\a(\b)=\gamma$. Remarquons que comme
$H^\pi_\beta(1)=H^{\tilde s_\a\pi}_{s_\a\b}(1)$,  la
  fonction $H^\pi_\beta$ est positive sur $[0,1]$. Donc pour $t\geq 1/2$,
  $H^{\tilde s_\a\pi}_{s_\a\b}(t)=H^\pi_\beta(t)$ est positif. 
  Si $t\leq 1/2$, alors~:
$$
H^{\tilde s_\alpha\pi}_{{s_\a\b}}(t)=H^\pi_{s_\alpha\beta}(t)
=H^\pi_\beta(t)-H^\pi_\alpha(t)\scal\alpha\bt.
$$ 
Si $\scal\alpha \bt \leq 0$, alors l'expression ci-dessus est
bien positive pour $t\leq 1/2$. 
Et si $\scal\alpha \bt > 0$, on vérifie que
$\scal{\pi(1)}{s_\alpha\bt}$ est positif et donc par hypothèse, pour tout
$t\in[0,1]$, $H^{\pi}_{s_\a\b}(t)\geq 0$.

Soit $\gamma$ une racine réelle et positive telle que
$H^{\tilde s_\a\pi}_{\g}(1)<0$. Alors il existe une racine réelle positive
$\beta$ telle que $s_\a\b=\g$. Comme
$H^{\pi}_{\b}(1)=H^{\tilde s_\a\pi}_{{s_\a\b}}(1)$, il existe $t^\pi_\b$ tel que
pour $t\geq t^\pi_\b$, $H^{\pi}_{\b}(t)$ est décroissante et
strictement négative.  
Si $H_{s_\a\b}^\pi(1)\geq 0$, ou si $H_{s_\a\b}^\pi(1)<0$ et
$t_{s_\a\b}^\pi\geq 1/2$, alors
$H_{s_\a\b}^{\tilde
  s_\a\pi}(1/2)=H^\pi_{s_\a\b}(1/2)=H^\pi_{\b}(1/2)\geq 0$ et donc 
$t^\pi_\b\geq 1/2$. Comme pour
$t\leq 1/2$, $H^{\tilde s_\a\pi}_{s_\a\b}(t)=H^\pi_{s_\a\b}(t)$ et pour
$t\geq 1/2$, $H^{\tilde s_\a\pi}_{s_\a\b}(t)= H^\pi_\b(t)$, on en
déduit que $H^{\tilde s_\a\pi}_{s_\a\b}(t)$ est positive ou nulle pour
$t\leq t^\pi_\b$ et décroissante et strictement négative pour $t>t^\pi_\b$.    
Si $t_{s_\a\b}^\pi <1/2$, on a alors~: 
$H_{s_\a\b}^{\tilde
  s_\a\pi}(1/2)=H_{s_\a\b}^{\pi}(1/2)=H^\pi_\b(1/2)<0$, et donc
$t^\pi_\b\leq 1/2$ et la fonction $H_{s_\a\b}^{\tilde s_\a\pi}(t)$ est
positive ou nulle pour $t\leq t_{s_\a\b}^\pi$ et strictement négative
et décroissante pour $t>t_{s_\a\b}^\pi$.

\end{preuve}

\vspace*{.5cm}

Lorsque $\pi$ est un chemin
affine par morceaux, le critère ci-dessus peut s'écrire sous la forme
pratique qui suit. 

\begin{Cor}\label{crit:pratique}
  Soit $\pi$ un chemin tel qu'il existe des poids
  $(\chi_1,\chi_2,\ldots,\chi_p)\in X_\Q^p$ tels que
  $\pi=\pi_{\chi_1}*\pi_{\chi_2}*\cdots*\pi_{\chi_p}$ ; si pour toute
  racine réelle positive $\beta$ on a~: 
$$\forall j\in\{1,\ldots,p-1\} \ {\rm tel\  que\ 
    }\scal{\chi_1+\chi_2+\cdots+\chi_j}\bt<0\Rightarrow
\scal{\chi_{j+1}}\bt\leq 0,$$
alors $\pi$ est extrémal.
\end{Cor}

\begin{NBs}
\begin{enumerate}
\item Il est facile de vérifier que le critère du théorème
 \ref{th:main} n'est pas une condition nécessaire pour qu'un chemin
 $\pi$ soit extrémal. Nous donnons un exemple en type $A_2$ dans la
 figure 2~: $\pi$ est extrémal puisque $\tilde s_{\a_1}\pi$ est un
 chemin dominant mais ne vérifie pas le critère du théorème
 \ref{th:main} ; nous utilisons dans cette figure les notations et
 conventions de \cite{Bou}. 
\item Il est également facile de vérifier que si $\pi$ est extrémal
  alors $\pi$ vérifie la
  condition suivante : pour toute racine réelle positive $\beta$,
  
 \begin{enumerate}
    \item si $H^\pi_\b(1)\geq
  0$ alors pour tout $t\in[0,1]$,  $H^\pi_\beta(t)\geq 0$ ;
\item  si $H^\pi_\b(1)<0$, alors pour tout $t\in[0,1]$,
  $H^\pi_\beta(t)\geq \scal{\pi(1)}\bt$. 
 \end{enumerate}
Mais cette condition n'est pas une condition suffisante, comme le
montre l'exemple dans la figure 3~: $\pi$ vérifie les conditions
ci-dessus, mais n'est pas extrémal puisque $\tilde s_{\a_2}\pi$ n'est
pas un chemin dominant, bien que le poids  $\tilde s_{\a_2}\pi(1)$
soit dominant.

\end{enumerate}
\end{NBs}

\section{Les composantes PRV généralisées}\label{part:3}
Nous allons maintenant énoncer le résultat principal de cet article.

\begin{Th}\label{th:prv}
  Soit $\lieg$ une algèbre de Kac-Moody symétrisable, soient
  $(\mu,\nu)$ un couple de poids dominants,  $(v,w)$ un couple
d'éléments du groupe de Weyl,  soit $\beta$ une racine
positive telle que l'une des conditions suivantes soit
satisfaite~:
\begin{enumerate}[(i)]
\item $\beta$ est une racine simple et $\lieg$ est de dimension finie ;
\item $v^{-1}\beta$ est une racine simple ; 
\item $w^{-1}\beta$ est une racine simple.
\end{enumerate}
Soit $k$ un entier vérifiant les inégalités~: 
$$0\leq k\leq\min\{\scal{v\mu}\bt,\scal{w\nu}\bt\},$$ 
alors si le poids
$\lambda=v\mu+w\nu-k\beta$ est dominant la représentation
irréductible $V(\lambda)$ est de multiplicité non nulle dans le
produit tensoriel $V(\mu)\otimes V(\nu)$. 
\end{Th}

\begin{preuve}
Nous reprenons ici la preuve donnée dans \cite{MPR2} qui montre
comment se ramener au cas où $v^{-1}\beta$ est
simple.
Supposons que le théorème soit vrai dans le cas
{\it (ii)} et soit $\mu,\nu,\lambda,v,w,k$ et $\beta$ {\it simple}
vérifiant les hypothèses du théorème avec $\lieg$ de dimension finie. 
On a alors~:
\begin{eqnarray}
\lambda=v\mu+w\nu-k\beta.\label{eq:linesym}
\end{eqnarray}
On peut transformer
cette égalité ( où $w_0$ désigne l'élément de plus grande longueur)~:
\begin{eqnarray}
\mu&=&v^{-1}\lambda+v^{-1}ww_0(-w_0\nu)+kv^{-1}\beta\\
&=&v^{-1}\lambda+s_{v^{-1}\b}v^{-1}ww_0(-w_0\nu)-(\scal{v^{-1}w\nu}{v^{-1}\bt}-k)v^{-1}\beta.\label{eq:linesym3}
\end{eqnarray}
Posons  $\lambda'=\mu$, $\mu'=\lambda$, $\nu'=-w_0\nu$, 
$v'=v^{-1}$, $w'=s_\alpha v^{-1}ww_0$, $\beta'=v^{-1}\beta$ et $k'=\scal{v^{-1}w\nu}{\alpha^\vee}-k$, l'équation (\ref{eq:linesym3})
devient~: 
\begin{eqnarray}
\lambda'=v'\mu'+w'\nu'-k'\alpha\,.\label{eq:linesym4}
\end{eqnarray}
On vérifie que les hypothèses sur $k$ se traduisent par~:
\begin{eqnarray*}
k'\geq 0,\\
k'\leq \scal{v'\nu'}{\bt}\\
k'\leq \scal{ w'\lambda'}{\bt}.
\end{eqnarray*}  
Et donc,  les hypothèses du théorème \ref{th:prv} sont satisfaites
dans le cas {\it (ii)} ($v'^{-1}\beta'$ est une racine simple) et
donc $V(\lambda')(\simeq V(\mu))$ est de multiplicité non nulle dans
$V(\mu')\otimes V(\nu')\simeq V(\lambda)\otimes V^*(\nu)$, ce qui
implique que $V(\lambda)$ est une composante de $V(\mu)\otimes V(\nu)$.
L'équivalence entre les cas {\it (ii)} et {\it (iii)} est immédiate en
utilisant la commutativité du produit tensoriel et sans hypothèse restrictive sur
$\lieg$.

On supposera maintenant que $\alpha=v^{-1}\beta$ est une racine simple. 
Nous allons construire un chemin de $B(\pi_\mu)*B(\pi_\nu)$ de poids
$$
v^{-1}\lambda=\mu+v^{-1}w\nu-k\a
$$ 
et extrémal ce qui impliquera bien le théorème, d'après
\ref{th:decompo:representation}. 
Soit $l:=\scal{w\nu}\bt=\scal{v^{-1}w\nu}\at$. 
Comme $0\leq k\leq l$ et d'après le
point 4 de la proposition \ref{prop:simple} le
chemin $f^{k}_\a\pi_{v^{-1}w\nu}$ est défini ; de plus si
$a=\frac{k}{l}$, on a alors~:  
$$
f^{k}_\a\pi_{v^{-1}w\nu}=\pi_{as_{\a}v^{-1}w\nu}*\pi_{(1-a)v^{-1}w\nu}.
$$
D'après la proposition \ref{prop:ls}, le chemin $\pi_\mu*f^{k}_\a\pi_{v^{-1}w\nu}$
appartient à $B(\pi_\mu)*B(\pi_\nu)$, et
pour montrer que ce chemin est
extrémal, nous allons utiliser le corollaire~\ref{crit:pratique}. Les
points de changements de
direction de $\pi_\mu*f^{k}_\a\pi_{v^{-1}w\nu}$ sont $\mu$ et
$\mu+as_{\a}v^{-1}w\nu$. 
Soit $\g$ une racine réelle positive ; remarquons d'abord que comme
$\mu$ est dominant, on ne peut pas avoir $\scal{\mu}\gt<0$. 
D'autre part, si $\g=\a$
alors $\scal{\mu+as_{\a}v^{-1}w\nu}{\a}\geq 0$, en effet~:
$$
\scal{\mu+as_{\a}v^{-1}w\nu}{\a}
=\scal\mu{\a^\vee}-a\scal{v^{-1}w\nu}\at=(\scal\mu{\a^\vee}-k).  
$$
Et le membre de droite est bien positif puisque
$\scal\mu\at=\scal{v\mu}\bt$. 
Supposons maintenant que $\g\neq \a$ et
$\scal{\mu+as_{\a}v^{-1}w\nu}\gt<0$, il faut montrer que
$a\scal{v^{-1}w\nu}{\gt}\leq 0$. 
Pour cela remarquons que~:
$$
\scal{\mu+as_{\a}v^{-1}w\nu}\gt=\scal{\mu}{\gt}-k\scal\a\gt+a\scal{v^{-1}w\nu}\gt
$$
Si $\scal\a\gt\leq 0$, alors~:
$$
\scal{\mu}{\gt}-k\scal\a\gt\geq 0
$$ 
et donc $a\scal{v^{-1}w\nu}{\gt}\leq 0$.
Si $\scal\a\gt> 0$, alors on a~:
\begin{eqnarray*}
\scal{\mu}{\gt}-k\scal\a\gt+a\scal{v^{-1}w\nu}{\gt}& & \\
 &\geq&\scal{\mu}\gt-\scal\mu\at\scal\a\gt+a\scal{v^{-1}w\nu}\gt\\
 &=&\scal{\mu}{s_\a\gt}+a\scal{v^{-1}w\nu}\gt.
\end{eqnarray*}
Comme $\g\neq\a$, $s_\a\g$ est une racine positive,
$\scal{\mu}{s_\a\g}$ est positif, et l'inégalité 
$a\scal{v^{-1}w\nu}{\gt}\leq 0$ est vérifiée.
\end{preuve}

\begin{NBs}
  \begin{enumerate}
\item Le cas $k=0$ correspond à l'énoncé PRV original. Dans ce cas le
  chemin à considérer est le chemin $\pi_\mu*\pi_{v^{-1}w\nu}$ qui est
  évidemment extrémal puisqu'il n'admet qu'un seul changement de
  direction en un poids dominant. C'est ce chemin que considère
  P. Littelmann dans sa preuve de la conjecture PRV dans \cite{pet1}. 
\item  Nous illustrons le théorème \ref{th:main} dans un cas ou
  $\lieg$ est de de 
    type $G_2$. Nous reprenons les notations de
    \cite{Bou}. En prenant $\mu=2\varpi_2$,
    $\nu=2(\varpi_1+\varpi_2)$, $v$ et $w$ de sorte que
    $v^{-1}w\nu=-8\varpi_1+2\varpi_2$, $p=1$, $\beta=3\a_1+\a_2$ et
    $k=1$, on obtient 
    $\lambda=\varpi_1+\varpi_2$ et donc $V(\varpi_1+\varpi_2)$ est une
    composante de 
    $V(\mu)\otimes V(\nu)$. Sur la figure 3, nous avons représenté le
    chemin extrémal 
    $\pi=\pi_\mu*f_{\a_2}\pi_{v^{-1}w\nu}$ ainsi que les chemins 
     extrémaux intermédiaires entre $\pi$ et
   le chemin dominant $\eta$ (en gras) tel que $\eta(1)=\lambda$.
  \end{enumerate}
 \end{NBs}

Nous concluons cet travail en énonçant un résultat plus général,
résultat énoncé et démontré dans \cite{MPR2} dans le cas fini. 

\begin{Th}\label{th:prvbis}
Soit $\lieg$ une algèbre de Kac-Moody symétrisable, soient
$\mu,\nu$ deux poids dominants de $\lieg$, $(v,w)\in W^2$, et
$(\beta_1,\beta_2,\ldots,\beta_p)$, $p$ racines 
orthogonales deux à deux. Supposons qu'une des trois assertions
suivantes soit vraie~:
\begin{enumerate}[(i)]
\item pour tout $i\in\{1,2,\ldots,p\}$, $\beta_i$ est une racine simple
  et $\lieg$ est de dimension finie ;
\item pour tout $i\in\{1,2,\ldots,p\}$, $v^{-1}\beta_i$ est une racine
  simple ;
\item pour tout $i\in\{1,2,\ldots,p\}$, $w^{-1}\beta_i$ est une racine
  simple.
\end{enumerate}
Soient $(k_1,k_2,\ldots,k_p)$ $p$ entiers et soit $\lambda$ le poids
$\lambda=v\mu+w\nu-\sum_{i=1}^p k_i\beta_i$. 
Supposons que pour pour $i\in\{1,2,\ldots,p\}$, les inégalités
suivantes soient vérifiées~:
$$
0\leq k_i\leq m_i=\min\{\scal{v\mu}\bti,\scal{w\nu}\bti\}
$$
alors si $\lambda$ est dominant, la représentation
$V(\lambda)$ est une composante de $V(\mu)\otimes V(\nu)$. 
\end{Th}

\begin{preuve}
La preuve de ce théorème étant une adaptation de la preuve du
théorème \ref{th:prv}, nous ne donnerons pas tous les détails.
En utilisant les mêmes arguments que dans la preuve du théorème
précédent, on peut supposer que pour tout $i\in\{1,2,\ldots,p\}$,
$\alpha_i=v^{-1}\beta_i$ est une racine simple.
Pour $i\in\{1,2,\ldots,p\}$ posons
$l_i:=\scal{w\nu}\bti=\scal{v^{-1}w\nu}\ati$ ; comme $0\leq k_i\leq
l_i$, d'après le
point 4 de la proposition \ref{prop:simple} le
chemin $f^{k_i}_\a\pi_{v^{-1}w\nu}$ est défini et de plus si
$a_i=\frac{k_i}{l_i}$, alors on a l'égalité~:  
$$
f^{k_i}_\a\pi_{v^{-1}w\nu}=\pi_{a_is_{\a_i}v^{-1}w\nu}*\pi_{(1-a_i)v^{-1}w\nu}.
$$
Comme les $(\beta_i)_{1\leq
  i\leq p}$ sont deux à deux orthogonales, il en est de même des
$(\alpha_i)_{1\leq i\leq p}$ et les
  opérateurs $f_{\alpha_i}$ commutent deux à deux. 
Le chemin
$$
\pi:=f^{k_1}_{\alpha_1}f^{k_2}_{\alpha_2}\ldots
f^{k_s}_{\alpha_p}\pi_{v^{-1}w\nu}
$$
est donc défini. Si on suppose que l'on
  a ordonné les $(\beta_i)_{1\leq i\leq p}$ de sorte que~: $
a_1\leq a_2\leq\cdots\leq a_p$,
alors le chemin $\pi$ s'écrit~:
$$
\pi_{a_1s_{\a_p}\ldots
  s_{\a_1}v^{-1}w\nu}*\pi_{(a_2-a_{1})s_{\a_p}\ldots
  s_{\a_2}v^{-1}w\nu}*\cdots *\pi_{(a_p-a_{p-1})s_{\a_p}v^{-1}w\nu}*\pi_{(1-a_p)v^{-1}w\nu}
$$
avec $a_0=0$ et $a_{p+1}=1$.
Les points de changement de direction pour
$\pi$ sont inclus dans l'ensemble
$\{\theta_1,\theta_2,\ldots,\theta_p\}$ avec 
\begin{eqnarray}
\theta_i&=&\sum_{j=1}^i(a_j-a_{j-1})s_{\a_p}\ldots
s_{\a_j}v^{-1}w\nu\nonumber\\
&=&\sum_{j=1}^i a_js_{\a_p}\ldots
s_{\a_j}v^{-1}w\nu-\sum_{j=1}^i a_{j-1}s_{\a_p}\ldots s_{\a_j}v^{-1}w\nu\nonumber\\ 
&=&a_is_{\a_p}\ldots s_{\a_i}v^{-1}w\nu+\sum_{j=1}^{i-1}a_j(s_{\a_p}\ldots s_{\a_j}v^{-1}w\nu-s_{\a_p}\ldots s_{\a_{j+1}}v^{-1}w\nu)\nonumber\\
&=&a_is_{\a_p}\ldots
s_{\a_i}v^{-1}w\nu-\sum_{j=1}^{i-1}a_j\scal{v^{-1}w\nu}{\at_j}\a_j\nonumber\\ 
&=&a_is_{\a_p}\ldots s_{\a_i}v^{-1}w\nu-\sum_{j=1}^{i-1}k_j\a_j\label{eq:1}\\ 
&=&a_i(v^{-1}w\nu-\sum_{j=i+1}^pl_j\a_j)-\sum_{j=1}^{i-1}k_j\a_j.\label{eq:2}
\end{eqnarray}

Le chemin $\pi_\mu*\pi$ est de poids $v^{-1}\lambda$ et il suffit de
montrer que ce chemin est 
extrémal ; pour cela nous allons utiliser le corollaire
\ref{crit:pratique}. Les
points de changements de
direction de $\pi_\mu*\pi$ sont inclus dans l'ensemble~:
$\{\mu,\mu+\theta_1,\ldots,\mu+\theta_p\}$.  
Soit $\g$ une racine réelle positive ; remarquons d'abord que comme
$\mu$ est dominant, on ne peut pas avoir $\scal{\mu}\gt<0$. 
D'autre part, s'il existe $m\in\{1,2,\ldots,p\}$ tel que $\g=\a_m$
alors $\scal{\mu+\theta_i}{\a_m}\geq 0$. En effet, en utilisant
l'égalité (\ref{eq:2}), si $m\leq i-1$ alors~:
$$
\scal{\mu+\theta_i}{\a_m}=\scal\mu{\a^\vee_m}+a_il_m-2k_m
=(\scal\mu{\a^\vee _m}-k_m)+(a_il_m-k_m).
$$
Les deux termes du membre de droite sont positifs, le  premier par
hypothèse et le
deuxième car $a_i\geq a_m$.
Si $m>i-1$, alors~:
$$
\scal{\mu+\theta_i}{\a_m}=\scal\mu{\a^\vee_m}+a_il_m-2a_il_m\geq\scal\mu{\a^\vee_m}-a_ml_m=\scal\mu{\a^\vee_m}-k_m
$$
et ce dernier terme est bien positif par hypothèse.
Supposons maintenant que $\g\not\in\{\a_1,\a_2,\ldots,\a_p\}$ et
que $\scal{\mu+\theta_i}\gt<0$, soit d'après l'égalité (\ref{eq:1}) 
$$
\scal{a_is_{\a_p}\ldots s_{\a_i}v^{-1}w\nu+\mu-\sum_{j=1}^{i-1}k_j\a_j}\gt<0.
$$
Il faut montrer que~: 
\begin{equation}\label{amontrer}
 \scal{s_{\a_p}\ldots s_{\a_{i+1}}v^{-1}w\nu}{\bt}\leq 0.  
\end{equation}
On peut supposer qu'il existe $i'\in\{0,1,\ldots,i\}$ tel que pour
$j\leq i'$, $\scal{\a_j}\gt\leq 0$ et pour $j>i'$,
$\scal{\a_j}\gt>0$. On a alors~:
\begin{eqnarray*}
 & & \scal{a_is_{\a_p}\ldots s_{\a_i}v^{-1}w\nu+\mu-\sum_{j=1}^{i-1}k_j\a_j}\gt\\
 &=&a_i\scal{s_{\a_p}\ldots s_{\a_{i+1}}v^{-1}w\nu}\gt
+\scal\mu\gt-\sum_{j=1}^{i}\scal{\a_j}\gt\\
&\geq&a_i\scal{s_{\a_p}\ldots s_{\a_{i+1}}v^{-1}w\nu}\gt
+\scal\mu\gt-\sum_{j=i'+1}^{i}\scal\mu\alpha_j^\vee\scal{\a_j}\gt\\
&=&a_i\scal{s_{\a_p}\ldots s_{\a_{i+1}}v^{-1}w\nu}\gt+\scal{\mu}{s_{\a_p}\ldots s_{\a_{i'+1}}\gt}.\\
\end{eqnarray*}
Comme $\g\not\in\{\a_1,\a_2,\ldots,\a_p\}$ et que les racines
$\a_1,\ldots \a_p$ sont orthogonales deux à deux, on a
$\scal{\mu}{s_{\a_p}\ldots s_{\a_{i'+1}}\gt}\geq 0$ et  
l'inéquation \ref{amontrer} est
  vérifiée.
\end{preuve}

\newpage

\bibliographystyle{amsalpha}
\bibliography{PRV_gen_chemin}

\newpage
\begin{figure}[ht!]\label{fig0}
 \centering
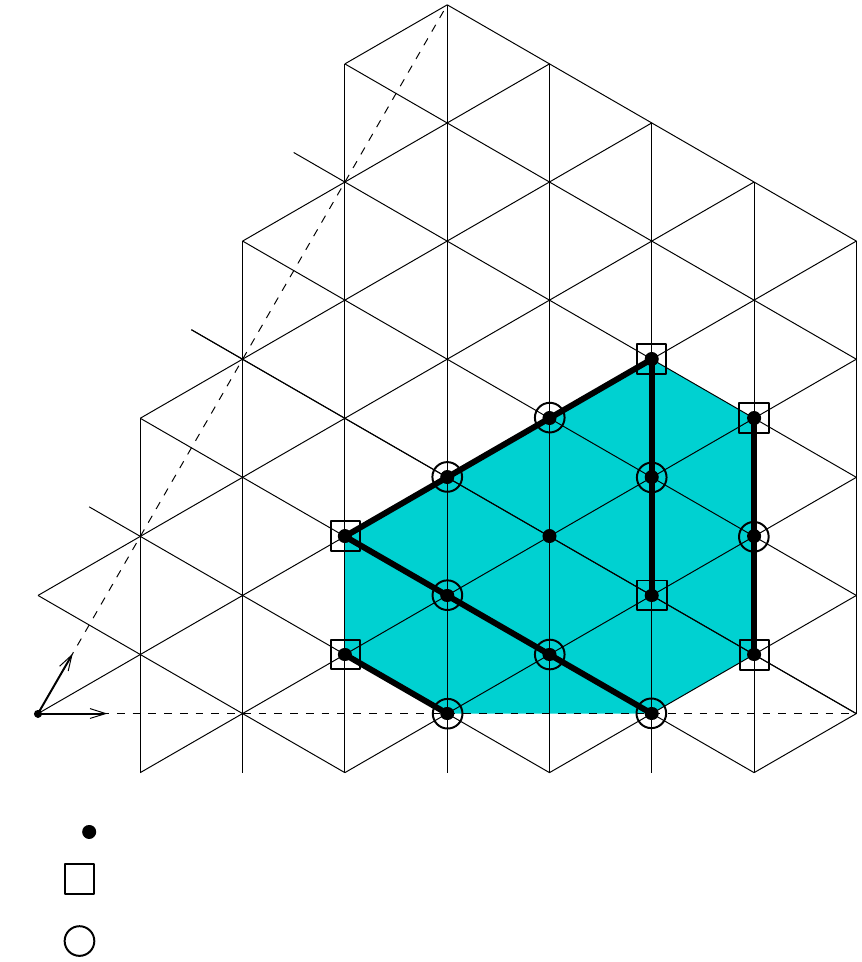 
\caption{}
\end{figure}

\begin{figure}[ht!] \label{fig1}
\centering
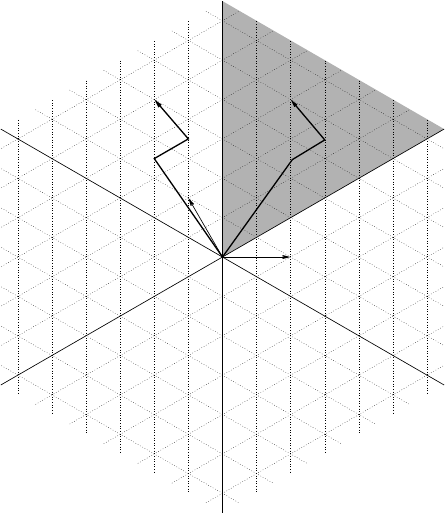  
\caption{}
\end{figure}  

\begin{figure}[ht!] \label{fig2}
\centering
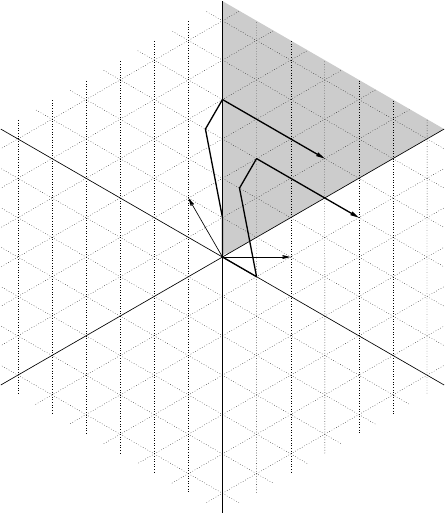 
\caption{} 
\end{figure}
\begin{figure}[ht!]\label{fig3}
 \centering
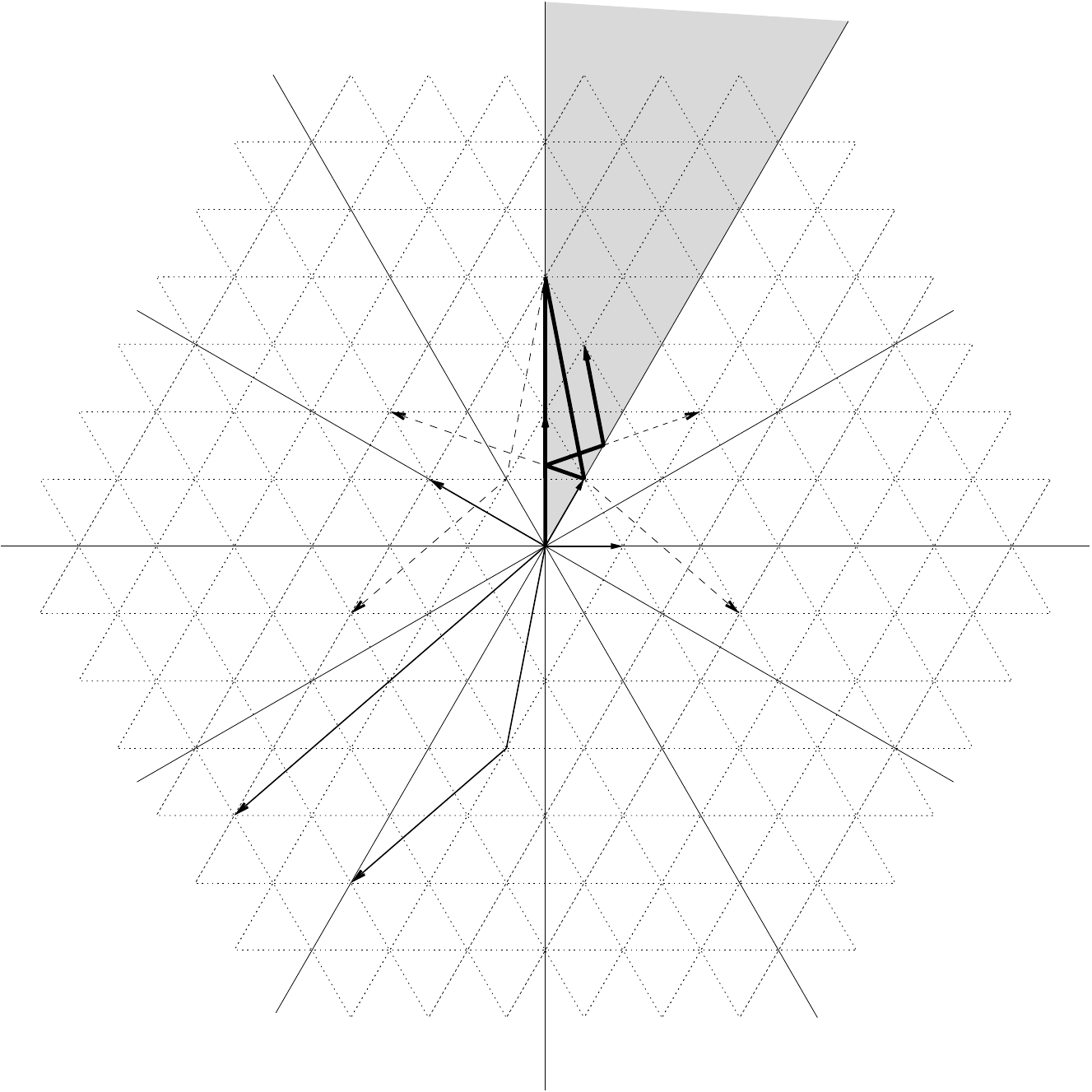 
\caption{}
\end{figure}

\end{document}